\documentclass[reqno,11pt]{amsart}
\usepackage{amsmath,amssymb,amsthm,graphicx,a4wide,enumerate,url}
\usepackage[small,bf]{caption}
\theoremstyle{plain}

\theoremstyle{definition}

\theoremstyle{remark}

\newcommand{\prn}[1]{\left(#1\right)}

\newcommand{\pd}[2]{\frac{\partial#1}{\partial#2}}
\newcommand{\ud}[1]{\,\mathrm{d}#1}
\newcommand{\xfoot}{\vec{x}_{\text{foot}}}
\newcommand{\tpdt}{t_{n+1}}
\newcommand{\tphdt}{t_{n+\frac{1}{2}}}
\begin{document}
\parskip.9ex

%===========================================================================
%=================================================================== Titles.
\title[A Comparative Study of the Efficiency of Jet Schemes]
{A Comparative Study of the Efficiency of Jet Schemes}
\author[P. Chidyagwai]{Prince Chidyagwai}
\address[Prince Chidyagwai]
{Department of Mathematics \\ Temple University \\ \newline
1805 North Broad Street \\ Philadelphia, PA 19122}
\email{chidyagp@temple.edu}
\urladdr{http://www.math.temple.edu/\~{}chidyagp}
\author[J.-C. Nave]{Jean-Christophe Nave}
\address[Jean-Christophe Nave]
{Department of Mathematics and Statistics \\ McGill University \\
\newline 805 Sherbrooke W. \\ Montreal, QC, H3A 2K6, Canada}
\email{jcnave@math.mcgill.ca}
\urladdr{http://www.math.mcgill.ca/jcnave}
\author[R. R. Rosales]{Rodolfo Ruben Rosales}
\address[Rodolfo Ruben Rosales]
{Department of Mathematics \\ Massachusetts Institute of Technology \\
\newline 77 Massachusetts Avenue \\ Cambridge, MA 02139}
\email{rrr@math.mit.edu}
\author[B. Seibold]{Benjamin Seibold}
\address[Benjamin Seibold]
{Department of Mathematics \\ Temple University \\ \newline
1805 North Broad Street \\ Philadelphia, PA 19122}
\email{seibold@temple.edu}
\urladdr{http://www.math.temple.edu/\~{}seibold}
%===========================================================================
\begin{abstract}
We present two versions of third order accurate jet schemes, which achieve high order accuracy by tracking derivative information of the solution along characteristic curves. For a benchmark linear advection problem, the efficiency of jet schemes is compared with WENO and Discontinuous Galerkin methods of the same order. Moreover, the performance of various schemes in tracking solution contours is investigated. It is demonstrated that jet schemes possess the simplicity and speed of WENO schemes, while showing several of the advantages as well as the accuracy of DG methods.
\end{abstract}
%===========================================================================

\maketitle

%===========================================================================
\section{Introduction}
\label{sec:introduction}
%===========================================================================
The advection of field quantities under a velocity field is an important sub-problem in many computational projects. Examples are the passive transport of concentrations, or the movement of interfaces using level set approaches \cite{OsherSethian1988}. We consider the linear advection equation
\begin{equation}
\label{eq:linear_advection}
\phi_t+\vec{v}\cdot\nabla\phi = 0\;,
\end{equation}
which moves a scalar quantity $\phi(\vec{x},t)$ by a given velocity field $\vec{v}(\vec{x},t)$. Equation \eqref{eq:linear_advection} is augmented with initial conditions $\phi(\vec{x},0) = \Phi(\vec{x})$ and boundary conditions. All data is assumed smooth in space and time. Furthermore, it is assumed that the problem at hand can be described by, or embedded into, a rectangular computational domain, equipped with a regular grid.

High order accurate numerical approximations of \eqref{eq:linear_advection} on a fixed grid are commonly based on schemes that employ polynomials of sufficient degree to interpolate smooth solutions with the required accuracy. One popular type of approach are finite difference ENO~\cite{ShuOsher1988} or WENO~\cite{LiuOsherChan1994} methods. These store approximations of the solution values at the grid points, and achieve a high order polynomial approximation by considering local neighborhoods that are several grid points wide. Another type of approach are discontinuous Galerkin (DG) \cite{ReedHill1973, CockburnShu1988, CockburnShu2001} methods. These achieve high order accuracy by storing a high degree polynomial approximation in each grid cell, and approximating the flux through cell boundaries based on a weak formulation of \eqref{eq:linear_advection}. Both WENO and DG methods are based on a semi-discretization of \eqref{eq:linear_advection}, and achieve high accuracy in time by using strong stability preserving (SSP) Runge-Kutta schemes \cite{ShuOsher1988, GottliebShu1998, GottliebShuTadmor2001}.

Both types of approaches incur problems and difficulties. Due to the wide stencils of WENO methods, it is challenging to preserve high order near boundaries. Furthermore, their non-locality poses difficulties for an effective parallelization, and their use in conjunction with adaptive grid approaches \cite{BergerOliger1984}. In contrast, in DG methods, communication is limited to neighboring cells. However, DG approaches are characterized by costly quadratures over grid cells and their edges, significant time step restrictions, and non-trivial implementation aspects.

A new class of approaches for \eqref{eq:linear_advection}, so called \emph{jet schemes}, has been proposed recently \cite{NaveRosalesSeibold2010, SeiboldRosalesNave2012}. The goal of these methods is to provide an attractive compromise between WENO and DG methods, possessing the optimal locality and good resolution properties of DG, while being close to the computational efficiency and ease of implementation of WENO schemes. In this paper we conduct a comparative study of the efficiency and accuracy of this new class of methods. We consider two versions of third order accurate jet schemes, and compare them with WENO and DG methods of the same order. The considered jet schemes are described in Sect.~\ref{sec:jet_schemes}. Then, in Sect.~\ref{sec:discontinuous_galerkin}, DG methods are outlined, and in Sect.~\ref{sec:weno}, information about the WENO schemes that we employ is given. Sect.~\ref{sec:numerical_results} shows the numerical results obtained when applying the three types of approaches to two versions of a benchmark problem. In Sect.~\ref{sec:discussion}, a discussion of the observations and theoretical estimates of the computational cost are given.

%===========================================================================
\section{Jet Schemes}
\label{sec:jet_schemes}
%===========================================================================
Jet schemes are based on an advect--and--project approach in function spaces. Given an approximation to the solution of \eqref{eq:linear_advection} at time $t_n$, an approximate solution at time $t_{n+1} = t_n+\Delta t$ is obtained by the time step
\begin{equation*}
\phi^{n+1} = P \circ A_{t_{n+1},t_n}\phi^n\;.
\end{equation*}
Here $A_{t_{n+1},t_n}$ is an approximate advection operator, defined by evolving the solution along characteristics using a numerical ODE solver, and $P$ is a projection operator, given by a piecewise Hermite interpolation based on parts of the jet of the solution at the points of a cartesian grid. The fact that the projection requires data only at grid points allows to use this approach as a numerical scheme: when evaluating the advection operator, only the characteristic curves that go through grid points at $t_{n+1}$ need to be considered. In addition, the evolution of derivatives of the solution along these characteristic must be found. As outlined below, the required spacial derivatives of the advection operator can be found by analytical differentiation (Sect.~\ref{subsec:analytical_differentiation}) or by approximations based on tracking multiple nearby characteristics (Sect.~\ref{subsec:eps-fd}).

Jet schemes can be constructed in any space dimension, and for any order of accuracy \cite{SeiboldRosalesNave2012}. For simplicity, here we only describe third order schemes in two space dimensions. We consider a rectangular computational domain $\Omega\subset\mathbb{R}^2$, equipped with a regular cartesian grid of grid size $h$.

%---------------------------------------------------------------------------
\subsection{Projection}
%---------------------------------------------------------------------------
In a grid cell $[a,a+h]\times [b,b+h]\subset\Omega$, let the vertices be indexed by a vector $\vec{q}\in\{0,1\}^2$, such that the vertex of index $\vec{q}$ is at position $\vec{x}_{\vec{q}} = (a+h\,q_1,b+h\,q_2)$. On each of the four vertices, let a vector of data be given by $\phi_{\vec{\alpha}}^{\vec{q}}\;\forall\,\vec{\alpha}\in\{0,1\}^2$. The data represents partial derivatives of orders up to 1 in each variable, as follows: $\vec{\alpha} = (0,0)$ represents function values $\phi$; $\vec{\alpha} = (1,0)$ and $\vec{\alpha} = (0,1)$ represent first derivatives $\partial_x\phi$ and $\partial_y\phi$, respectively; and $\vec{\alpha} = (1,1)$ represents $\partial_{xy}\phi$. This data is interpolated by the bi-cubic polynomial
\begin{equation}
\label{eq:bi-cubic_interpolant}
\mathcal{H}(\vec{x}) = \sum_{\vec{q},\vec{\alpha}\in\{0,1\}^p}
\phi_{\vec{\alpha}}^{\vec{q}}\;W_{\vec{\alpha}}^{\vec{q}}(\vec{x})\;,
\end{equation}
where the $W_{\vec{\alpha}}^{\vec{q}}(\vec{x})$ are bi-cubic basis functions, given by the tensor product formulas
\begin{equation*}
W_{\vec{\alpha}}^{\vec{q}}(\vec{x}) = h^{\alpha_1+\alpha_2}\;
w_{\alpha_1}^{q_1}\big(\tfrac{x-a}{h}\big)\,
w_{\alpha_2}^{q_2}\big(\tfrac{y-b}{h}\big)\;,
\end{equation*}
and the $w_{\alpha}^{q}$ are the univariate basis functions
\begin{align*}
w_0^0(x) = 1-3x^2+2x^3, \quad
w_0^1(x) =   3x^2-2x^3, \quad
w_1^0(x) = x-2x^2+ x^3, \quad
w_1^1(x) =   -x^2+ x^3.
\end{align*}
The bi-cubic interpolant \eqref{eq:bi-cubic_interpolant} is an $O(h^4)$ accurate approximation to any sufficiently smooth function $\phi$ that it interpolates on a cell of size $h$ \cite{SeiboldRosalesNave2012}.

On the computational domain $\Omega\subset\mathbb{R}^2$, let the grid points be labeled $\vec{x}_{\vec{m}}$, where $\vec{m}\in\mathbb{Z}^2$. For any (sufficiently smooth) function $\phi:\Omega\to\mathbb{R}$, we define a global interpolant $\mathcal{H}_{\phi}$ as follows: at each grid point $\vec{x}_{\vec{m}}$, evaluate the derivatives of $\phi$, $\partial_x\phi$, $\partial_y\phi$, and $\partial_{xy}\phi$, to produce a data vector $\phi_{\vec{\alpha}}^{\vec{m}}\;\forall\,\vec{\alpha}\in\{0,1\}^2$. Then, on each grid cell, use this data to define the bi-cubic interpolant \eqref{eq:bi-cubic_interpolant}. In the function space
\begin{equation*}
S^{2,+} = \{\psi\in C^1: \psi\text{~twice differentiable a.e. with~}D^2\psi \in L^\infty\}\;,
\end{equation*}
this procedure can be applied using the following convention: whenever $\partial_{xy}\psi$ must be evaluated at a point at which $D^2\psi$ is not defined in the classical sense, we define it as
\begin{equation*}
\partial_{xy}\psi(\vec{x}_{\vec{m}}) = \tfrac{1}{2}\Big(
\text{ess}\limsup_{\vec{x} \to \vec{x}_{\vec{m}}}\partial_{xy}\psi(\vec{x})+
\text{ess}\liminf_{\vec{x} \to \vec{x}_{\vec{m}}}\partial_{xy}\psi(\vec{x})\Big)\;.
\end{equation*}
The view on the general order of approximation case \cite{SeiboldRosalesNave2012} reveals the rationale for the notation $S^{2,+}$. Clearly, the re-application of the interpolation procedure does not change the result: $\mathcal{H}_{\mathcal{H}_{\phi}} = \mathcal{H}_{\phi}$. Hence, in the space $S^{2,+}$ it can be formulated as a projection operator
\begin{equation}
\label{eq:projection}
P\phi = \mathcal{H}_{\phi}\;.
\end{equation}

%---------------------------------------------------------------------------
\subsection{Advection}
%---------------------------------------------------------------------------
The characteristic form of equation \eqref{eq:linear_advection} is
\begin{equation}
\label{eq:characteristic_system}
\frac{d\phi}{dt} = 0
\quad\text{along}\quad
\frac{d\vec{x}}{dt} = \vec{v}(\vec{x},t)\;.
\end{equation}
Let $\vec{X}(\vec{x},\tau,t)$ denote the solution of the ODE for the characteristic curves at time $t$, when starting with initial conditions $\vec{x}$ at time $\tau$, i.e.~it is defined by
\begin{equation*}
\pd{}{t}\vec{X}(\vec{x},\tau,t) = \vec{v}(\vec{X}(\vec{x},\tau,t),t)
\quad\text{with}\quad
\vec{X}(\vec{x},\tau,\tau) = \vec{x}\;.
\end{equation*}
Then due to \eqref{eq:characteristic_system}, the solution of \eqref{eq:linear_advection} satisfies $\phi(\vec{x},\tau) = \phi(\vec{X}(\vec{x},\tau,t),t)$. In practice, the characteristic ODE \eqref{eq:characteristic_system} must be approximated. Let $\vec{\mathcal{X}}(\vec{x},\tau,t)$ represent an approximate solution of the ODE for the characteristic curves at time $t$, when starting with initial conditions $\vec{x}$ at time $\tau$. It typically arises from a numerical ODE solver, e.g.~a high order Runge-Kutta step. Furthermore, introduce the associated approximate advection operator $A_{\tau,t}$, which maps the solution at time $t$ to an approximate solution at time $\tau$. It acts on a function $g(\vec{x})$ as follows
\begin{equation*}
(A_{\tau,t}\,g)(\vec{x}) = g(\vec{\mathcal{X}}(\vec{x},\tau,t))\;.
\end{equation*}
As shown in \cite{NaveRosalesSeibold2010}, the use of a locally $k^\text{th}$ order Runge-Kutta scheme results in a $k^\text{th}$ order accurate approximation to the solution
\begin{equation*}
A_{t+\Delta t,t}\,\phi(\vec{x},t)-
\phi(\vec{X}(\vec{x},t+\Delta t,t),t) = O(|\Delta t|^k)\;.
\end{equation*}
This means that through the characteristic equations \eqref{eq:characteristic_system}, any ODE solver induces an approximate advection operator of the same order of accuracy. However, this idea alone cannot be used as a numerical time stepping scheme, since the step from $t$ to $t+\Delta t$ in general generates a function $A_{t+\Delta t,t}\,g$ that cannot be represented with a finite amount of data. Therefore, at the end of every time step, we apply the projection operator \eqref{eq:projection}, which generates a function that can be stored on a computer. The function space $S^{2,+}$ is invariant under diffeomorphisms, thus both $P$ and $A_{t+\Delta t,t}$ map from $S^{2,+}$ into itself. Consequently, one full approximate solution step is given by applying $P\circ A_{t+\Delta t,t}$ to the approximate solution at time $t$.

%---------------------------------------------------------------------------
\subsection{Derivative Updates by Analytical Differentiation}
\label{subsec:analytical_differentiation}
%---------------------------------------------------------------------------
The application of the projection \eqref{eq:projection} requires the knowledge of $\psi$, $\partial_x\psi$, $\partial_y\psi$, $\partial_{xy}\psi$ at grid points, with $\psi = A_{t+\Delta t,t}\,\phi^n$. These spacial derivatives of $A_{t+\Delta t,t}\,\phi^n$ can be found by analytically differentiating the ODE solver, as demonstrated below for the Shu-Osher scheme \cite{ShuOsher1988}. One step with this scheme is $O((\Delta t)^4)$ accurate, which for the scaling $\Delta t\propto h$ matches the $O(h^4)$ accuracy of the projection operator \eqref{eq:projection}. Using the notation $t_n = t$, $\tpdt = t+\Delta t$, and $\tphdt = t+\frac{1}{2}\Delta t$, we obtain the updates for the required parts of the jet $\phi$, $\nabla\phi = (\partial_x\phi,\partial_y\phi)$, and $\partial_{xy}\phi$ at a grid point $\vec{x}$ as:
\begin{equation}
\label{eq:shu_osher}
\hspace{-.5em}
\begin{split}
\vec{x}_1 &= \vec{x} - \Delta t\,\vec{v}(\vec{x},\tpdt) \\
\nabla\vec{x}_1 &= I - \Delta t\,\nabla\vec{v}(\vec{x},\tpdt) \\
\partial_{xy}\vec{x}_1 &= -\Delta t\,\partial_{xy}\vec{v}(\vec{x},\tpdt) \\
\vec{x}_2 &= \tfrac{3}{4}\vec{x} + \tfrac{1}{4}\vec{x}_1 -
\tfrac{1}{4}\Delta t\,\vec{v}(\vec{x}_1,t_n) \\[-.1em]
\nabla\vec{x}_2 &= \tfrac{3}{4}I + \tfrac{1}{4}\nabla\vec{x}_1 -
\tfrac{1}{4}\Delta t\,\nabla\vec{x}_1\cdot\nabla\vec{v}(\vec{x}_1,t_n) \\[-.1em]
\partial_{xy}\vec{x}_2 &= \tfrac{1}{4}\partial_{xy}\vec{x}_1 -
\tfrac{1}{4}\Delta t ( \partial_{xy}\vec{x}_1\cdot\nabla\vec{v}(\vec{x}_1,t_n)
+((\partial_x\vec{x}_1)^T\!\cdot (\partial_y\vec{x}_1)):D^2\,\vec{v}(\vec{x}_1,t_n) ) \\
\xfoot &= \tfrac{1}{3}\vec{x} + \tfrac{2}{3}\vec{x}_2 -
\tfrac{2}{3}\Delta t\,\vec{v}(\vec{x}_2,\tphdt) \\[-.1em]
\nabla\xfoot &= \tfrac{1}{3}I + \tfrac{2}{3}\nabla\vec{x}_2 -
\tfrac{2}{3}\Delta t\,\nabla\vec{x}_2 \cdot
\nabla\vec{v}(\vec{x}_2,\tphdt) \\[-.1em]
\partial_{xy}\xfoot &= \tfrac{2}{3}\partial_{xy}\vec{x}_2 -
\tfrac{2}{3}\Delta t (
\partial_{xy}\vec{x}_2\cdot\nabla\vec{v}(\vec{x}_2,\tphdt)
+((\partial_x\vec{x}_2)^T\!\cdot (\partial_y\vec{x}_2)):D^2\,
\vec{v}(\vec{x}_2,\tphdt) ) \hspace{-1.5em} \\
\phi(\vec{x},\tpdt) &= \mathcal{H}(\xfoot,t_n) \\
(\nabla\phi)(\vec{x},\tpdt) &= \nabla\xfoot\cdot\nabla\mathcal{H}(\xfoot,t_n) \\
(\partial_{xy}\phi)(\vec{x},\tpdt) &=
\partial_{xy}\xfoot\cdot\nabla\mathcal{H}(\xfoot,t)+
((\partial_x\xfoot)^T\!\cdot (\partial_y\xfoot)):D^2\mathcal{H}(\xfoot,t_n)
\end{split}
\end{equation}
In this approach, the characteristic curve is tracked from $\vec{x}$ at time $t+\Delta t$ back to $\xfoot$ at time $t$. The data at this position is given by the Hermite interpolation, defined by the data (at time $t$) at the four vertices of the cell that $\xfoot$ is contained in. The update rules for the derivatives are systematically inherited from the update rule for the function value, and they match exactly what one would obtain when evolving the solution using $A_{t+\Delta t,t}$ everywhere, and then applying the spacial derivatives.

%---------------------------------------------------------------------------
\subsection{Order of Accuracy and Stability}
%---------------------------------------------------------------------------
One step of the jet scheme described above is fourth order accurate, since the advection operator is $O((\Delta t)^4)$ accurate, and the projection operator is $O(h^4)$ accurate. As usual, when going from this local error to the global error, one order is lost, since $O(\tfrac{1}{\Delta t})$ time steps are required to reach the final time. Hence, the presented scheme is globally third order, given that it is stable.

Like in many other numerical approaches, stability is not automatically guaranteed. As shown in \cite{SeiboldRosalesNave2012}, jet schemes can be constructed (by using a different projection than \eqref{eq:projection}) that are unstable. However, the jet scheme based on the projection \eqref{eq:projection} is stable. A key factor in the stability argument is that among all (sufficiently smooth) functions that match given data on a cartesian grid, the Hermite interpolant \eqref{eq:bi-cubic_interpolant} minimizes the stability functional
\begin{equation*}
\mathcal{F}[\phi] = \int_\Omega (\partial_{xxyy}\phi(\vec{x}))^2\ud{\vec{x}}\;.
\end{equation*}
Hence $\mathcal{F}[P\phi] \le \mathcal{F}[\phi]$ for all sufficiently smooth $\phi$, which yields bounds on the amounts of oscillations that the numerical scheme can create \cite{SeiboldRosalesNave2012}.

%---------------------------------------------------------------------------
\subsection{A Simpler and More Efficient Derivative Tracking}
\label{subsec:eps-fd}
%---------------------------------------------------------------------------
The analytical differentiation procedure, outlined with example \eqref{eq:shu_osher}, is the most general approach to derive the update rules for derivatives. However, it is not always the simplest to implement. An alternative approach, which is solely based on tracking function values, is provided by \emph{$\varepsilon$-finite differences}. For the here considered third order jet scheme in 2D, the following procedure can be employed. As shall be seen in Sect.~\ref{sec:numerical_results}, it leads to a very simple and efficient numerical scheme.

At a grid point $\vec{x} = (x,y)$, instead of tracking one characteristic curve from time $t_{n+1}$ back to time $t_n$, we track four characteristic curves, starting at $\vec{x}^{\vec{q}} = (x+q_1\varepsilon,y+q_2\varepsilon)$ where $\vec{q}\in\{-1,1\}^2$, and $\varepsilon$ is a small number, see below. Let the corresponding characteristic foot-points be denoted $\xfoot^{\vec{q}}$. The updates for the required derivatives are obtained as follows:
\begin{enumerate}[(1)]
\item
The ``center of mass'' $\tfrac{1}{4}\big(\xfoot^{(1,1)}+\xfoot^{(-1,1)}+\xfoot^{(1,-1)}+\xfoot^{(-1,-1)}\big)$ determines to which cell all four foot-points are associated to.
\item
The Hermite interpolant \eqref{eq:bi-cubic_interpolant} that corresponds to that cell is evaluated at each foot-point to yield the values $\phi(\vec{x}^{\vec{q}},t_{n+1}) = \mathcal{H}(\xfoot^{\vec{q}},t_n)$.
\item
The required parts of the jet are defined as
\begin{align*}
\phi(\vec{x},t_{n+1}) &=
\;\,\tfrac{1}{4}\,(\phi(\vec{x}^{(1,1)},t_{n+1})+\phi(\vec{x}^{(-1,1)},t_{n+1})
+\phi(\vec{x}^{(1,-1)},t_{n+1})+\phi(\vec{x}^{(-1,-1)},t_{n+1})) \\
\partial_x\phi(\vec{x},t_{n+1}) &=
\,\tfrac{1}{4\varepsilon}\,(\phi(\vec{x}^{(1,1)},t_{n+1})-\phi(\vec{x}^{(-1,1)},t_{n+1})
+\phi(\vec{x}^{(1,-1)},t_{n+1})-\phi(\vec{x}^{(-1,-1)},t_{n+1})) \\
\partial_y\phi(\vec{x},t_{n+1}) &=
\,\tfrac{1}{4\varepsilon}\,(\phi(\vec{x}^{(1,1)},t_{n+1})+\phi(\vec{x}^{(-1,1)},t_{n+1})
-\phi(\vec{x}^{(1,-1)},t_{n+1})-\phi(\vec{x}^{(-1,-1)},t_{n+1})) \\
\partial_{xy}\phi(\vec{x},t_{n+1}) &=
\tfrac{1}{4\varepsilon^2}(\phi(\vec{x}^{(1,1)},t_{n+1})-\phi(\vec{x}^{(-1,1)},t_{n+1})
-\phi(\vec{x}^{(1,-1)},t_{n+1})+\phi(\vec{x}^{(-1,-1)},t_{n+1}))
\end{align*}
\end{enumerate}
All approximations are $O(\varepsilon^2)$ accurate. In addition, round-off errors of up to $O(\delta/\varepsilon^2)$ arise, where $\delta$ is the accuracy of the floating point operations. With the optimal choice $\varepsilon = O(\delta^{1/4})$, the errors incurred by the $\varepsilon$-finite differences are of magnitude $O(\delta^{1/2})$. Thus, with double precision arithmetics, the presented approach can be used up to a desired accuracy of $10^{-7}$.

%===========================================================================
\section{Discontinuous Galerkin Method}
\label{sec:discontinuous_galerkin}
%===========================================================================
DG methods \cite{ReedHill1973, CockburnShu1988} have a wide area of application. They can successfully approximate nonlinear problems on unstructured geometries, and thus are much more general than the problems considered in this paper. However, since the task at hand is the numerical approximation of the linear advection equation \eqref{eq:linear_advection} on a regular grid, it is natural to investigate the accuracy and efficiency of DG methods for it. The DG methodology that we use here is in line with the ideas presented in \cite{Cockburn1997, CockburnShu2001}. For convenience, we restrict the presentation here to incompressible velocity fields $\nabla\cdot\vec{v} = 0$. In this case, the advection equation \eqref{eq:linear_advection} can be written in conservative form as
\begin{equation}
\label{eq:linear_advection_conservative}
\phi_t+\nabla\cdot(\vec{v}\phi) = 0\;.
\end{equation}
Equation \eqref{eq:linear_advection_conservative} is discretized in space as follows. Let $\mathcal{T}_h$ be a regular triangulation of the computational domain $\Omega$. At each instance in time, we seek an approximation $\phi_h$ of $\phi$, such that $\phi_h(t)$ belongs to the finite dimensional space
\begin{equation*}
S^k_h = \{\psi\in L^1(\Omega): \psi|_{K} \in P^k_h(K)\;\forall\,K\in\mathcal{T}_h\}\;,
\end{equation*}
where $P^k_h(K)$ denotes the space of polynomials of degree $\le k$ that live on the element $K$. Here, we choose $k=2$ to obtain a third order scheme. As in the standard DG formulation \cite{Cockburn1997}, we replace $\phi$ by $\phi_h$ in \eqref{eq:linear_advection_conservative}, multiply by a test function $\psi_h \in S^k_h(K)$, and integrate over $K\in\mathcal{T}_h$ to obtain the semi-discrete formulation
\begin{equation}
\label{eq:weak_form_upwind}
\frac{d}{dt}\int_{K}\phi_h\psi_h\ud{x}
= \int_{K}\vec{v}\cdot\nabla\psi_h\phi_h\ud{x}
-\int_{e\in\partial K} \widehat{\phi_h}\,\vec{v}\cdot\vec{n}_{e,K}\psi_h\ud{\Gamma}
\quad\forall\,\psi_h\in S^k_h(K)\;.
\end{equation}
Here $\vec{n}_{e,K}$ is the outward unit normal on edge $e$, and $\widehat{\phi_h}$ is an expression that comprises the values of $\phi_h$ on each side of the edge. We use an upwind flux, which is defined as
\begin{equation*}
\widehat{\phi_h} =
\begin{cases}
\phi_h^{+} & \text{if~} \vec{v}\cdot\vec{n}_{e,K} \ge 0 \\
\phi_h^{-} & \text{otherwise}\;,
\end{cases}
\end{equation*}
where $\phi_h^\pm = \lim_{\varepsilon\to 0^\pm}\phi_h(\vec{x}+\varepsilon\vec{n}_{e,K},t)$, and the velocity field is evaluated in the edge center. The integrals over edges and elements in \eqref{eq:weak_form_upwind} are approximated by Gaussian quadrature rules that are exact for polynomials of degree $2k$ over the elements, and exact over polynomials of degree $2k+1$ over the edges. This results in an ODE system of the form
\begin{equation}
\label{eq:dg_ode_system}
\frac{d}{dt}\phi_h = L_h(\phi_h)
\quad\text{with}\quad
\phi_h(\vec{x},0) = \Phi_h(\vec{x})\;,
\end{equation}
where $L_h(\phi_h)$ is the spacial discretization of the operator $-\nabla\cdot(\vec{v}\phi)$, and $\Phi_h(\vec{x})\in S^k_h$ approximates $\Phi(\vec{x})$. The time integration of \eqref{eq:dg_ode_system} is done using the third order SSP Shu-Osher scheme \cite{ShuOsher1988}. Stability is ensured if the specific CFL condition
\begin{equation*}
\Delta t < \tfrac{1}{c(2k+1)}\,h
\end{equation*}
with $c>1$ is guaranteed \cite{CockburnShu2001}. For the examples considered in this paper, the choice $c=2$ turns out to yield the lowest cost vs.~accuracy ratio. Since for the numerical tests conducted here the solutions are very smooth, no slope limiters are implemented.

The DG code used here is based on triangular elements. In order to have a fair comparison with the cartesian grids that WENO and jet schemes use, we design meshes as follows. For a desired resolution $h$, first a cartesian grid of resolution $\sqrt{2}h$ is constructed. Then each cell is divided along its diagonal into two triangles. Like a cartesian grid of resolution $h$, the resulting triangular mesh consists of $1/h^2$ elements, and the shortest height in an element is $h$.

%===========================================================================
\section{WENO Scheme}
\label{sec:weno}
%===========================================================================
Like DG methods, WENO schemes \cite{LiuOsherChan1994} are based on a semi-discretization of \eqref{eq:linear_advection} in space, and the use of SSP schemes \cite{ShuOsher1988, GottliebShu1998} to advance in time. Unlike DG methods, WENO schemes are specific to regular grids. We employ the third order WENO finite difference approach described in \cite{JiangShu1996}, using the Shu-Osher scheme \cite{ShuOsher1988} with $\Delta t = h$ in time. For the spacial approximation, equation \eqref{eq:linear_advection} is rewritten as
\begin{equation*}
\phi_t = -u\,\partial_x\phi-v\,\partial_y\phi\;.
\end{equation*}
At each grid point, both derivatives $\partial_x\phi$ and $\partial_y\phi$ are approximated in a univariate fashion by using values on four grid points in each coordinate direction: two in the upwind direction, and one in the downwind direction. WENO schemes have limiters built into the derivative approximations: smoothness indicators yield a nonlinear weighted average of different polynomial approximations of the derivatives. We consider two versions of limiting: one as described in \cite{JiangShu1996}, and another one without limiting (which yields a linear finite difference scheme).

%===========================================================================
\section{Numerical Results}
\label{sec:numerical_results}
%===========================================================================
We test the accuracy, the relative efficiency, and the quality of tracking contours of jet schemes, DG methods, and WENO approaches using the classical \emph{vortex in a box} flow test \cite{BellColellaGlaz1989, LeVeque1996}, adapted as follows. On the computational domain $(x,y)\in [0,1]^2$, and for $t\in [0,T]$, we consider the linear advection equation \eqref{eq:linear_advection} with the velocity field
\begin{equation}
\label{eq:swirl}
\vec{v}(x,y,t) = \cos\prn{\pi\tfrac{t}{T}}
\begin{pmatrix}
\phantom{-}\sin^2(\pi x)\,\sin(2\pi y) \\
-          \sin(2\pi x) \,\sin^2(\pi y)
\end{pmatrix},
\end{equation}
which is a model for the passive swirling and successive un-swirling of a concentration field by an incompressible fluid motion. This test is a mathematical analog of well-known ``unmixing'' experiments \cite{Heller1960, FluidMixingMovie}. For this flow, the final solution equals the initial conditions, i.e.~$\phi(x,y,T) = \phi(x,y,0)$. We consider smooth initial conditions, and periodic boundary conditions. Note that for linear advection problems, this test is quite general. Moreover, it has a built-in difficulty parameter: the larger $T$, the more challenging it is to numerically resolve the highly elongated contours of the solution.

\begin{figure}
\begin{minipage}[b]{.49\textwidth}
\includegraphics[width=\textwidth]{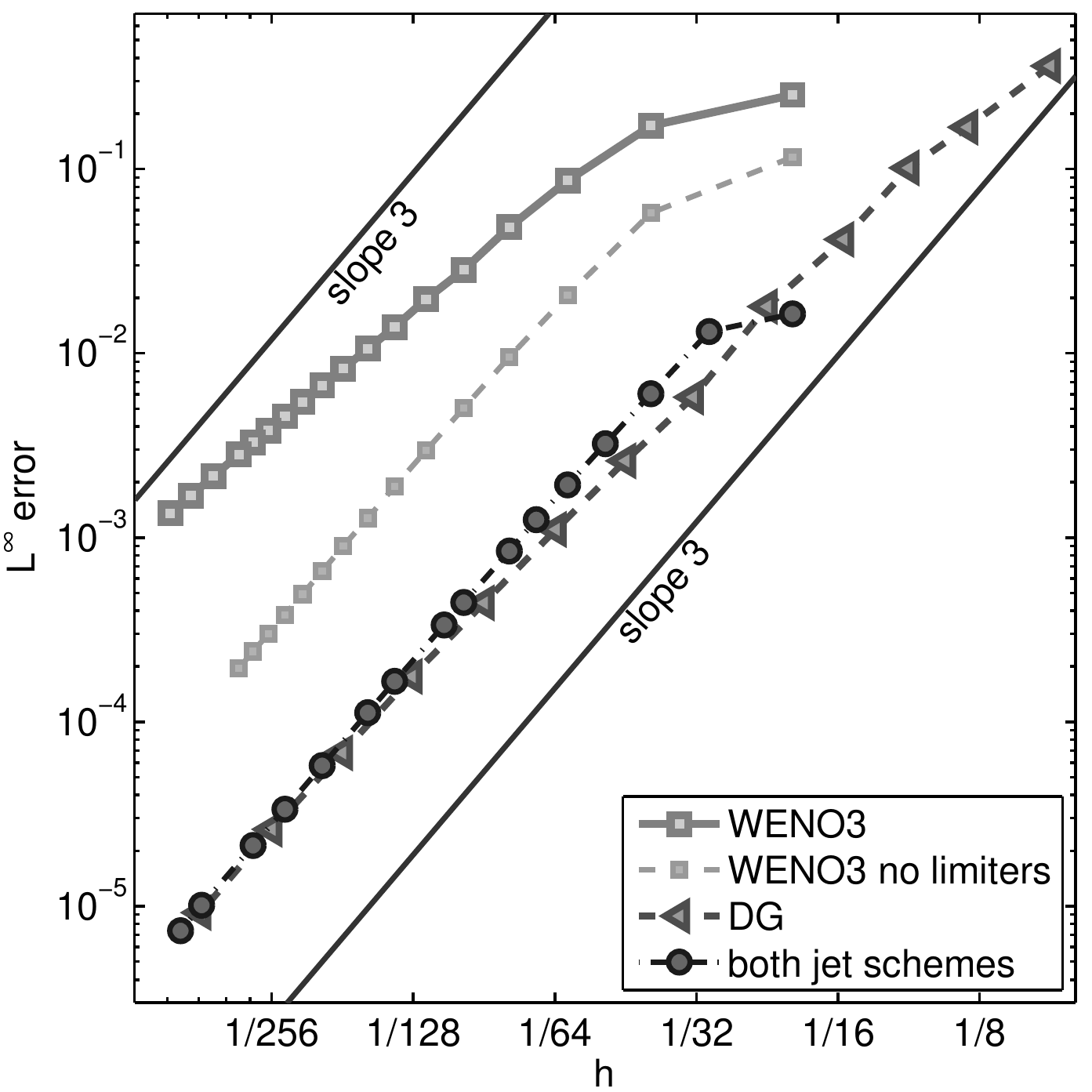}
\vspace{-1.8em}
\caption{Error convergence for third order jet schemes, DG, and WENO.}
\label{fig:comparison_convergence}
\end{minipage}
\hfill
\begin{minipage}[b]{.49\textwidth}
\includegraphics[width=\textwidth]{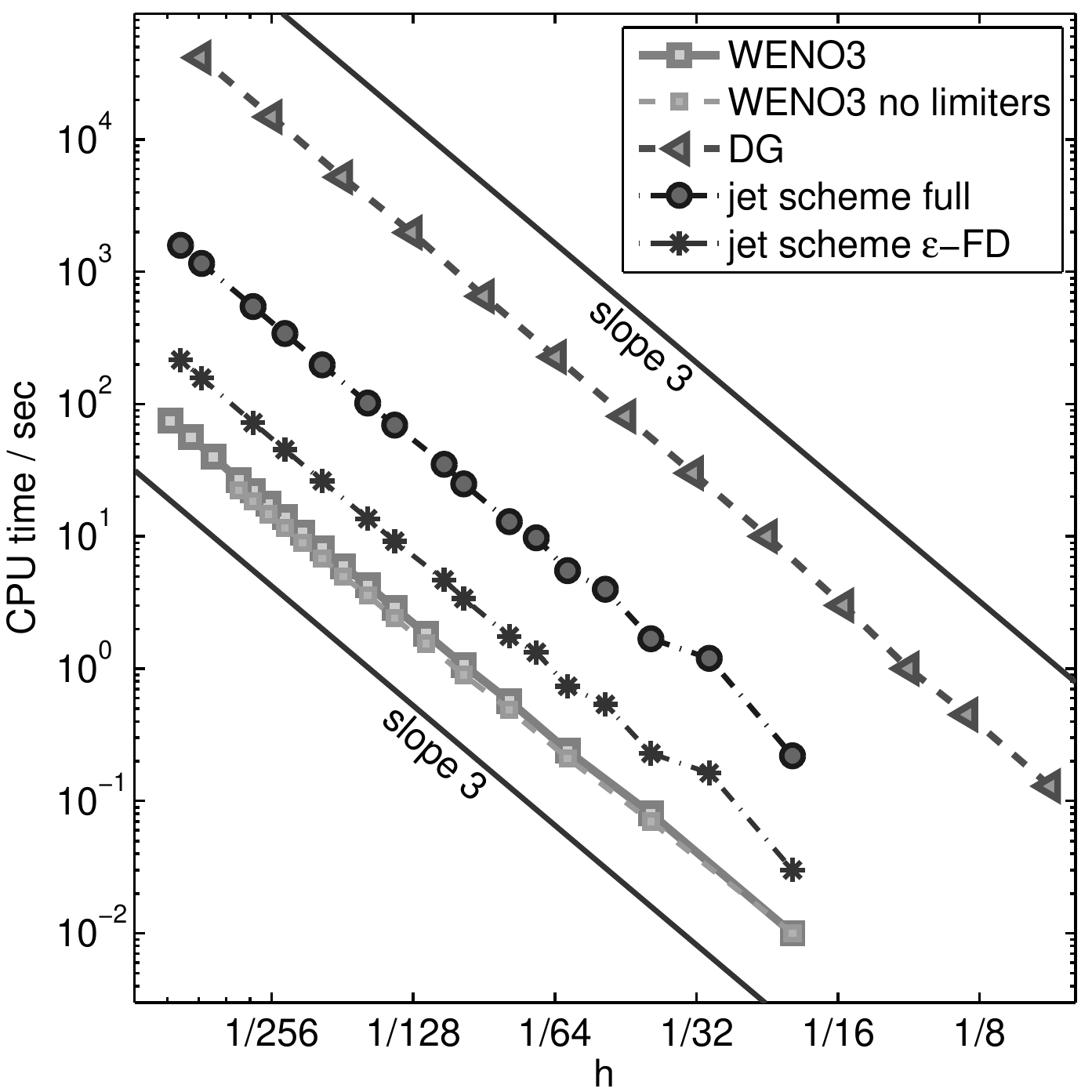}
\vspace{-1.8em}
\caption{Scaling of the computational cost for jet schemes, DG, and WENO.}
\label{fig:comparison_cost}
\end{minipage}

\vspace{1em}
\begin{minipage}[b]{.49\textwidth}
\includegraphics[width=\textwidth]{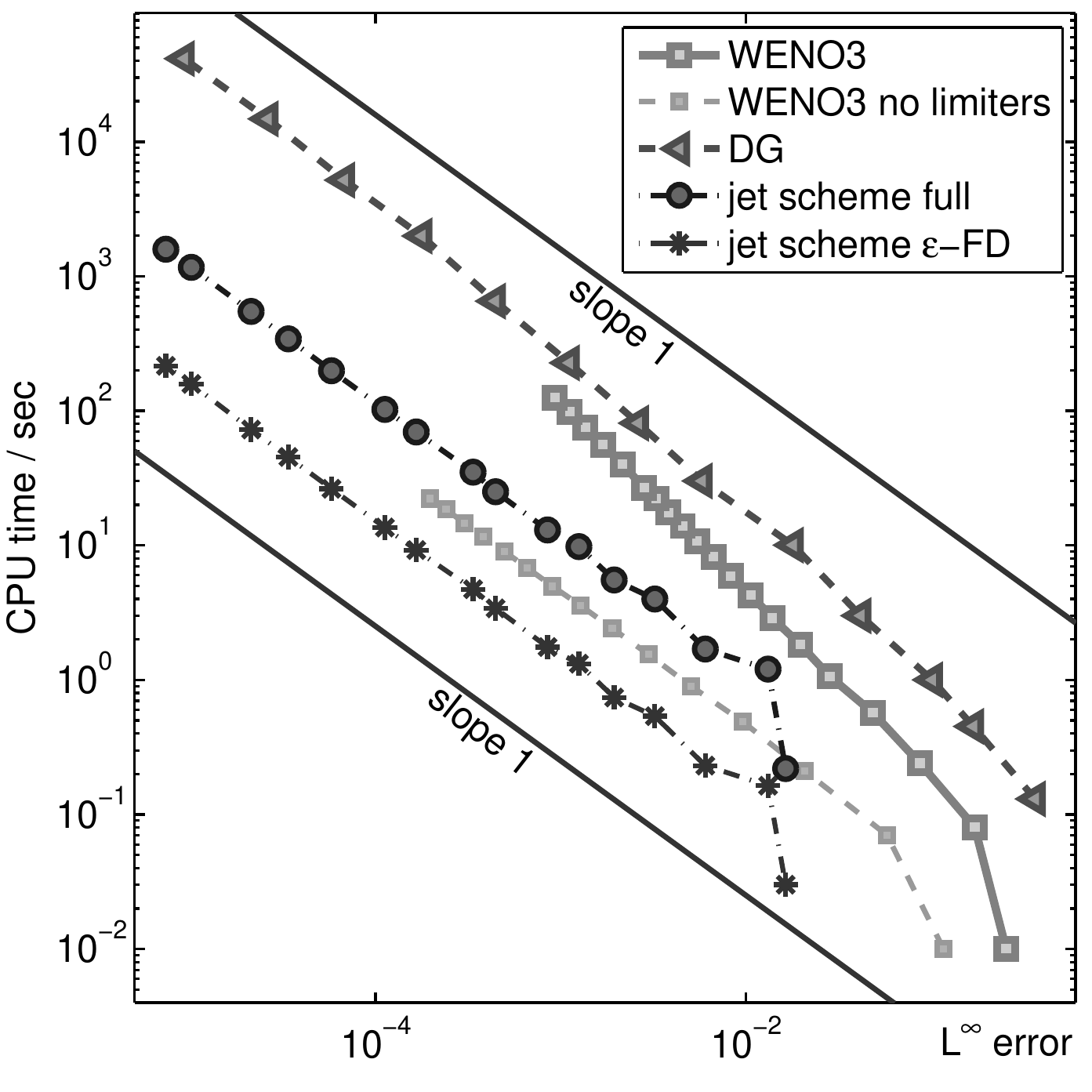}
\vspace{-1.8em}
\caption{Efficiency of jet schemes, DG, and WENO, measured as cost over accuracy.}
\label{fig:comparison_efficiency}
\end{minipage}
\end{figure}

%---------------------------------------------------------------------------
\subsection{Numerical Efficiency}
\label{subsec:efficiency}
%---------------------------------------------------------------------------
For the investigation of numerical efficiency, we consider smooth periodic initial conditions $\phi(x,y,0) = \cos(2\pi x)\cos(4\pi y)$, and choose $T = 1$. For various mesh resolutions $h$, we apply five different numerical schemes to the test problem. The results are shown in Figs.~\ref{fig:comparison_convergence}--\ref{fig:comparison_efficiency}. Specifically, we consider:
\begin{enumerate}[(i)]
\item A classical third order WENO scheme (see Sect.~\ref{sec:weno}), denoted by large squares;
\item A third order WENO scheme without limiting (see Sect.~\ref{sec:weno}), denoted by small squares;
\item A third order DG method (see Sect.~\ref{sec:discontinuous_galerkin}), denoted by triangles;
\item A third order jet scheme with a full derivative update (see Sect.~\ref{subsec:analytical_differentiation}), denoted by circles;
\item A third order jet scheme based on $\varepsilon$-finite differences (see Sect.~\ref{subsec:eps-fd}), denoted by stars.
\end{enumerate}
For all methods, we measure the error at $t = T$ in the $L^\infty$ norm, and the CPU time. Figure~\ref{fig:comparison_convergence} shows the error as a function of the resolution $h$. One can observe that jet schemes and DG have roughly the same accuracy. In contrast, for the same resolution, WENO schemes yield significantly larger errors. It is also visible that the WENO3 with limiters does not achieve the full third order accuracy on the range of resolutions under consideration (see \cite{JiangShu1996}). It is apparent that for the smooth solution at hand, limiters are not beneficial. Figure~\ref{fig:comparison_cost} shows the CPU time as a function of the resolution $h$. One can see a clear ranking in computational costs: WENO schemes are very fast, jet schemes are more costly, and DG is even more costly. It is visible that $\varepsilon$-finite differences are not only simpler to implement, but are also significantly faster than the jet scheme based on the full derivative update. Finally, Fig.~\ref{fig:comparison_efficiency} shows the CPU time as a function of the error. This cost vs.~accuracy ratio measures the true efficiency of the numerical schemes. The results show that the low computational cost of WENO schemes renders them preferable over DG. However, WENO does not possess the optimal locality that DG has. Jet schemes provide an interesting alternative. While the full update version is slightly less efficient than the non-limited WENO, the $\varepsilon$-finite differences version is even more efficient.

\begin{figure}
\begin{minipage}[b]{.49\textwidth}
\includegraphics[width=\textwidth]{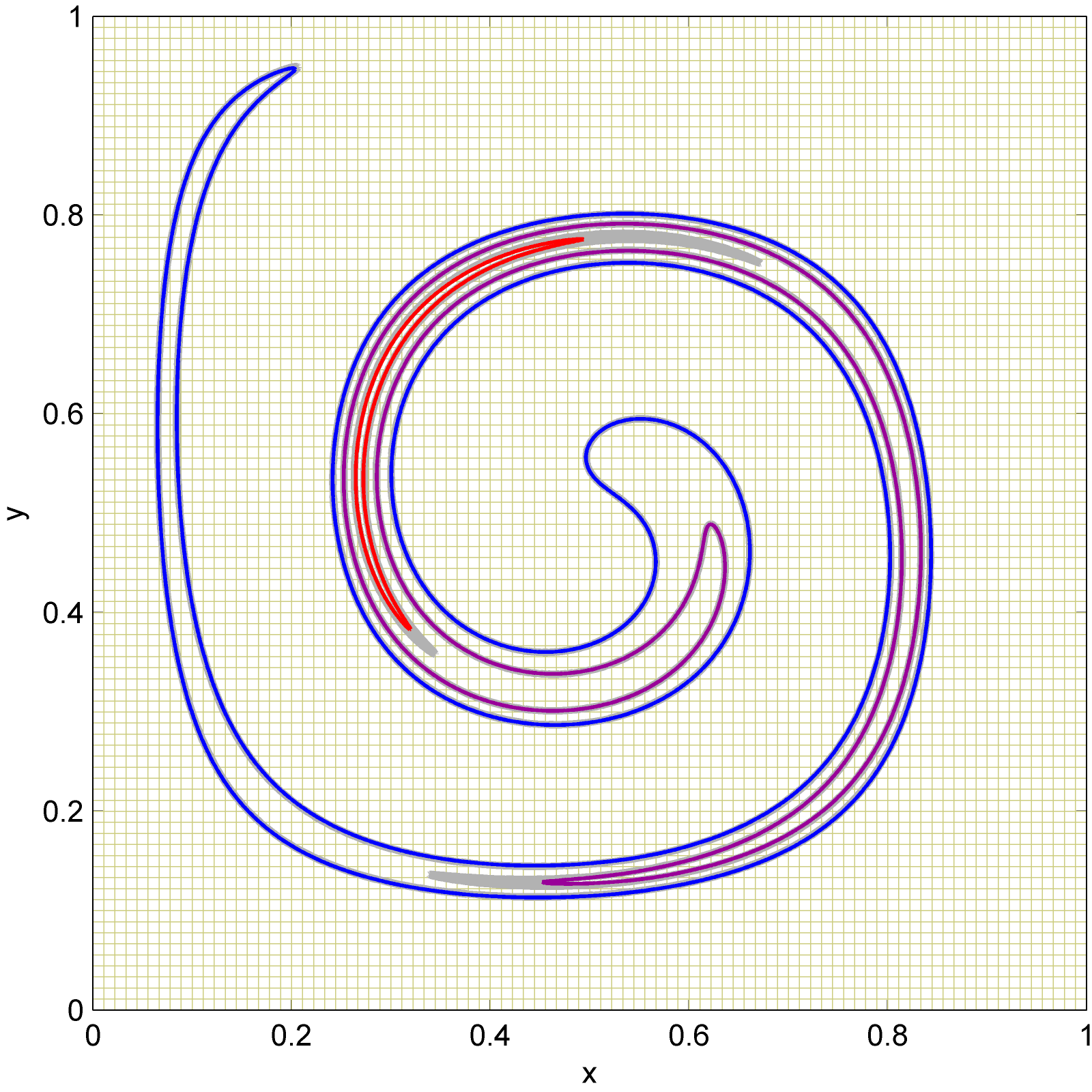}
\vspace{-1.8em}
\caption{Deformed contours, computed with a jet scheme using $h=1/90$.}
\label{fig:contours_full_jet}
\end{minipage}
\hfill
\begin{minipage}[b]{.49\textwidth}
\includegraphics[width=\textwidth]{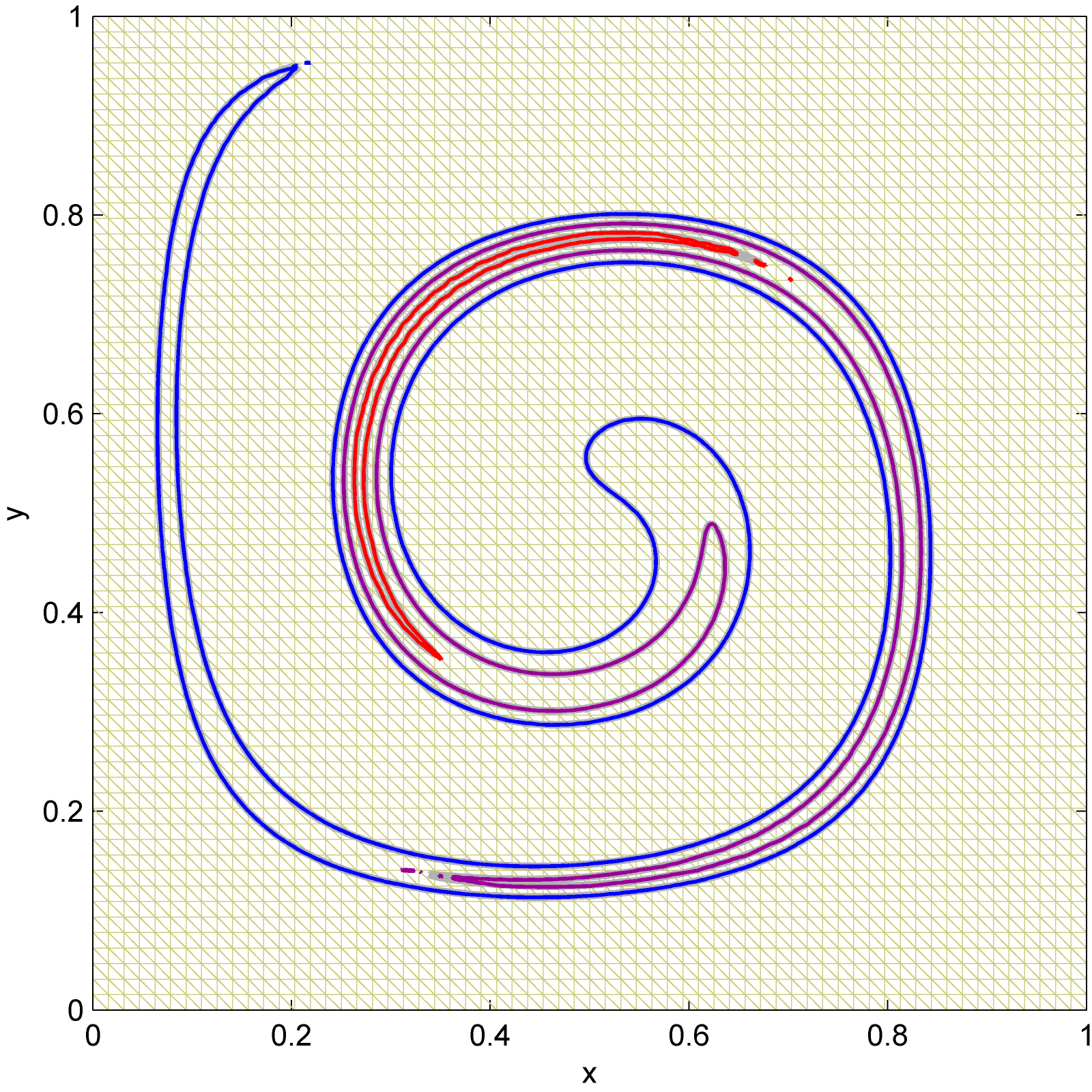}
\vspace{-1.8em}
\caption{Deformed contours, computed with a DG scheme using $h=1/90$.}
\label{fig:contours_full_dg}
\end{minipage}

\vspace{1em}
\begin{minipage}[b]{.49\textwidth}
\includegraphics[width=\textwidth]{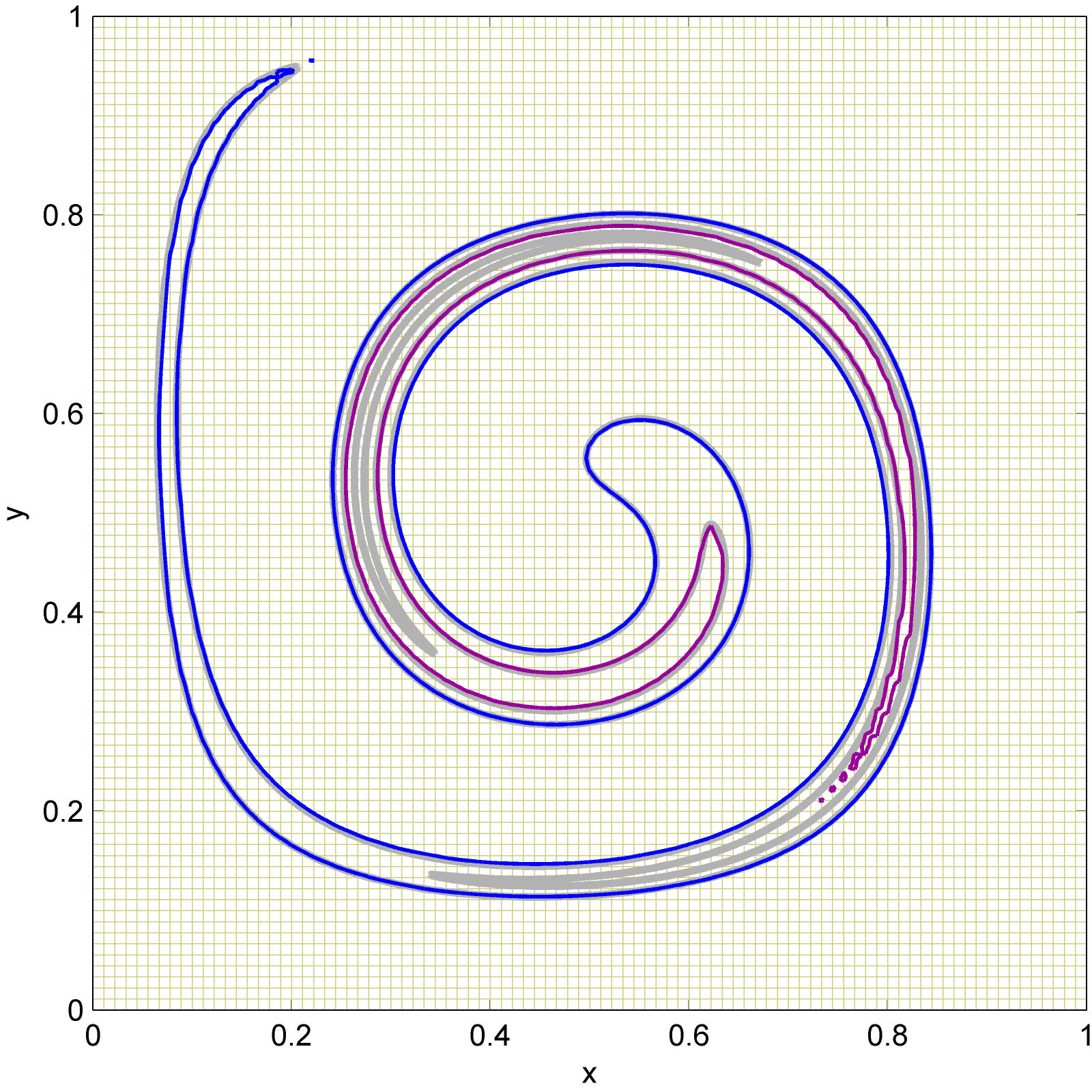}
\vspace{-1.8em}
\caption{Deformed contours, computed with a WENO scheme using $h=1/90$.}
\label{fig:contours_full_weno}
\end{minipage}
\end{figure}

%---------------------------------------------------------------------------
\subsection{Quality of Contours}
%---------------------------------------------------------------------------
In order to assess the quality of the numerical approximations that the different methods produce, we consider how accurately contours of the solution are deformed under the flow \eqref{eq:swirl} with $T = 6$. For this test, we choose the initial conditions
\begin{equation*}
\phi(x,y,0) = \exp\prn{-10(x-0.5)^2-10(x-0.75)^2}\;,
\end{equation*}
and consider three contours at $\phi = \exp(-10\,r_k^2)$, where $r_1 = 0.044$, $r_2 = 0.132$, and $r_3 = 0.220$. %$r\in\{0.044,0.132,0.220\}$.
At $t = 0$, these contours are three concentric circles of radii $r_k$, centered at the point $(0.5,0.75)$. We plot the deformed contours at $t = \frac{T}{2} = 3$, the time of their maximum deformation.

We compare a jet scheme, DG, and WENO for the same resolution $h$. As found in Sect.~\ref{subsec:efficiency}, these methods involve very different computational costs. Yet, this equal-$h$ comparison is of practical relevance, since frequently the advection equation \eqref{eq:linear_advection} is only a sub-problem in a larger project, and the computational grid is determined by the underlying application. Fig.~\ref{fig:contours_full_jet} shows the contours obtained with a jet scheme based on $\varepsilon$-finite differences (see Sect.~\ref{subsec:eps-fd}); Fig.~\ref{fig:contours_full_dg} shows the contours obtained with DG (see Sect.~\ref{sec:discontinuous_galerkin}); and Fig.~\ref{fig:contours_full_weno} shows the contours obtained with WENO without limiting (see Sect.~\ref{sec:weno}). In all cases, the inner (red) contour is $\phi = \exp(-10\,r_1^2)$, the middle (purple) contour is $\phi = \exp(-10\,r_2^2)$, and the outer (blue) contour is $\phi = \exp(-10\,r_3^2)$. Moreover, each contour underlays the contour of the true solution as a thick gray curve. In addition, the computational grid is shown in the background. Note that in all cases, contours are plotted using a sub-grid: for jet schemes, the solution inside a grid cell is given by the Hermite bi-cubic interpolant \eqref{eq:bi-cubic_interpolant}; for DG, the solution inside an element is given by a bivariate quadratic polynomial; and for WENO, we use a simple bi-linear interpolant on each grid cell.

As it was already visible in the error convergence graph shown in Fig.~\ref{fig:comparison_convergence}, also in terms of the tracking of contours, the jet scheme as well as DG are significantly more accurate than WENO. The comparison between jet schemes and DG is more interesting. Clearly, the DG contours shown in Fig.~\ref{fig:contours_full_dg} are closer to the true solution's contours than the jet scheme's contours are. However, at the tips of contours, DG produces small ripples and break-ups. In contrast, the jet scheme contours are extremely uniform, all the way to their end. Moreover, for the tracking of contours, DG has another potential disadvantage: due to the discontinuities of the numerical solution across element boundaries, the resulting contours may have small jumps (in fact, the contours in Fig.~\ref{fig:contours_full_dg} possess jumps, but they are too small to be visible to the eye). In certain applications, such as level set methods \cite{OsherSethian1988}, this could lead to problems. These observations come in addition to the fact that for the same resolution, DG is significantly more costly than jet schemes (see Fig.~\ref{fig:comparison_cost}).

%===========================================================================
\section{Discussion}
\label{sec:discussion}
%===========================================================================
The two types of traditional numerical methods, WENO and DG, achieve high order accuracy in very different ways: WENO considers function values in a wide neighborhood. In contrast, DG uses high order polynomials in each element. This gives DG a certain level of sub-grid resolution, and restricts communication to neighboring elements only. While based on a very different methodology, the recently proposed \cite{NaveRosalesSeibold2010, SeiboldRosalesNave2012} jet schemes share these advantages with DG.

The results presented here show that the advantages of DG come at the expense of efficiency. This observation can be explained by the following theoretical estimates. With a third order Runge-Kutta scheme, all presented schemes require three right hand side evaluations per time step. On a mesh of $N$ elements, even a fully optimized and problem-specific DG code requires at least $11.5N$ evaluations of the velocity field per right hand side evaluation (7 quadrature points on each of $N$ elements, and 3 quadrature points on each of $\frac{3}{2}N$ edges). In comparison, WENO required $N$ velocity field evaluations, and a jet scheme (using $\varepsilon$-finite differences) requires $4N$ velocity field evaluations (4 characteristics per grid point). In addition, due to its more restrictive CFL condition, DG requires 5--10 times more time steps than WENO or jet schemes. The results shown in Fig.~\ref{fig:comparison_cost} are in line with these theoretical estimates, albeit at more extreme ratios between the costs of the methods. These increased ratios are most likely due to additional data structure management requirements for DG, and also for jet schemes. However, the key results in the efficiency plot shown in Fig.~\ref{fig:comparison_efficiency} remain valid even if the measured CPU times were replaced by the theoretical cost estimates given above.

In terms of accuracy, jet schemes and DG yield similar accuracies (for the same number of elements). In contrast, the low cost of WENO is outweighed by its lower accuracy. Therefore, in terms of true efficiency, jet schemes and WENO are in the same range, while DG methods are significantly more costly. This observation is further confirmed by a comparison of the quality of contour approximation.

These investigations lead to the conclusion that jet schemes indeed fill a need as computational approaches for advection problems on regular grids, namely: they are simpler to implement, and less costly than DG. However, in contrast to WENO, they are optimally local and thus are preferable for the purpose of parallelization, adaptive mesh refinement, and the implementation of complicated boundary conditions. In addition, jet schemes impose no restrictions on the ODE solvers that are needed. In contrast, DG and WENO methods require SSP ODE solvers to ensure TVD stability \cite{GottliebShuTadmor2001}. The key weakness of jet schemes is their limited range of applicability. While established for linear advection \cite{NaveRosalesSeibold2010, SeiboldRosalesNave2012}, their application to more general advection problems is a the topic of current research.

%===========================================================================
\section*{Acknowledgments}
%===========================================================================
The authors would like to acknowledge the support by the National Science Foundation. J.-C. Nave, R. R. Rosales, and B. Seibold were supported through grant DMS--0813648, and P. Chidyagwai was supported through grant DMS--1115269. In addition, P. Chidyagwai, R. R. Rosales, and B. Seibold wish to acknowledge partial support by the National Science Foundation through grants DMS--1007967, DMS--1007899, DMS--1115269, and DMS--1115278. Further, J.-C. Nave wishes to acknowledge partial support by the NSERC Discovery Program. This research was supported in part by the National Science Foundation through major research instrumentation grant number CNS--09--58854.

%===========================================================================
\bibliographystyle{plain}
\bibliography{references_complete}
%===========================================================================

\end{document}